\documentclass[12pt,reqno]{amsart}
\begin{document}

\title[Injective envelopes and local multiplier algebras]{Injective envelopes
        and local multiplier algebras of C*-algebras}
\author[M.~Frank]{{\rm Michael Frank\\Universit\"at Leipzig, Mathematisches
Institut\\ D-04109 Leipzig, Fed.~Rep.~Germany}}
\email{frank@mathematik.uni-leipzig.de}

%%%%%%%%%%%%%%%%%%%%%%%%%%%%%%%%%%%%%%%%%%%%%%%%%%%%%%%%%%%%%%%%%%%%%%%%%%%%%
\begin{abstract}
The local multiplier C*-algebra $M_{loc}(A)$ of any C*-algebra $A$ can be
$*$-isomorphicly embedded into the injective envelope $I(A)$ of $A$ in
such a way that the canonical embeddings of $A$ into both these C*-algebras
are identified. If $A$ is commutative then $M_{loc}(A) \equiv I(A)$. The
injective envelopes of $A$ and $M_{loc}(A)$ always coincide, and every higher
order local multiplier C*-algebra of $A$ is contained in the regular monotone
completion $\overline{A} \subseteq I(A)$ of $A$.
For C*-algebras $A$ with a center $Z(A)$ such that $Z(A) \circ A$ is
norm-dense in $A$ the center of the local multiplier C*-algebra of $A$ is the
local multiplier C*-algebra of the center of $A$, and both they are
$*$-isomorphic to the injective envelope of the center of $A$.
A Wittstock type extension theorem for completely bounded bimodule maps on
operator bimodules taking values in $M_{loc}(A)$ is proven to hold if and
only if $M_{loc}(A) \equiv I(A)$. In general, a solution of the problem for
which C*-algebras $A$ the C*-algebras $M_{loc}(A)$ is injective
is shown to be equivalent to the solution of I.~Kaplansky's 1951 problem
whether all AW*-algebras are monotone complete.
\end{abstract}
%%%%%%%%%%%%%%%%%%%%%%%%%%%%%%%%%%%%%%%%%%%%%%%%%%%%%%%%%%%%%%%%%%%%%%%%%%%%%
\maketitle

\thispagestyle{empty}

The injective envelope of a C*-algebra in the category of C*-algebras and
completely positive linear maps is defined by an extrinsic algebraic
characterization.
M.~Hamana showed its general existence and uniqueness up to $*$-isomorphism,
cf.~\cite{Ham}. The main problem is to determine the injective envelope $I(A)$
of a given C*-algebra $A$ from the structure of $A$, i.e.~intrinsicly.
For commutative C*-algebras an intrinsic characterization was given by
H.~Gonshor in \cite{Go1,Go2}. He relied on I.~M.~Gel'fand's theorem for
commutative C*-algebras and on the topology of the locally compact Hausdorff
space $X$ corresponding to a commutative C*-algebra $A=C_0(X)$. Unfortunately,
there seems to be no obvious way to extend his results to the non-commutative
case. The injective envelope of non-commutative C*-algebras has been described
only for some examples and special classes of C*-algebras yet.

\smallskip
One of the goals of this short note is an intrinsic algebraic characterization
of the injective envelopes of commutative C*-algebras identifying them with
their local multiplier algebras. In the present paper we emphasize the
topological approach which relies on strict topologies. For the more algebraic
point of view see \cite{FP}. In the general non-commutative case we
show the existence of a canonical injective $*$-homomorphism mapping the local
multiplier algebra into the corresponding injective envelope (Th.~1).
The question whether the image of this map coincides with the entire injective
envelope has a negative answer, in general, and is connected with
I.~Kaplansky's still unsolved problem on the monotone completeness of
AW*-algebras in the case of AW*-algebras. So this class of C*-algebras cannot
be characterized in full at present. We indicate some examples (Cor.~7).
The center of the local multiplier C*-algebra of $A$ always coincides
with the local multiplier C*-algebra of the center $Z(A)$ of $A$, and both they can
be identified with the injective envelope of the center of $A$, in case $Z(A)
\circ A$ is dense in $A$, i.e.~$Z(M_{loc}(A)) \equiv M_{loc}(Z(A)) \equiv
I(Z(A))$ in that situation (Th.~2).
Moreover, the injective envelopes of $A$ and $M_{loc}(A)$ are always
$*$-isomorphic, and every higher order local multiplier C*-algebra of $A$ is
contained in the regular monotone completion $\overline{A} \subseteq I(A)$ of
$A$ (Th.~5).
Additionally, we establish extension properties of completely bounded maps
from operator spaces or operator (bi)modules into the (local) multiplier
C*-algebra of a certain C*-algebra (Prop.~4).

%%%%%%%%%%%%%%%%%%%%%%%%%%%%%%%%%%%%%%%%%%%%%%%%%%%%%%%%%%%%%%%%%%%%%%%%%%%%%%%

\section{Preliminaries}

A C*-algebra $B$ is said to be {\it injective} if given any self-adjoint
linear subspace $S$ of a unital C*-algebra $C$ with $1_C \in S$ any completely
positive linear map of $S$ into $B$ extends to a completely positive linear
map of $C$ into $B$, cf.~\cite[Def.~2.1]{Ham}. Equivalently, $B$ is injective
if and only if it is an injective object of the category consisting of
C*-algebras as objects and completely positive linear maps as morphisms.
For any C*-algebra $A$ we can define an {\it injective envelope} in the
following way (cf.~\cite{Ham}): let $B \supseteq A$ be an injective C*-algebra
such that $A$ is embedded into $B$ as a C*-subalgebra. If the identity map on
$A$ possesses a unique extension to a completely positive linear map on $B$
then $B$ is said to be an injective envelope $I(A)$ of $A$. The injective
envelope $I(A)$ of a C*-algebra $A$ is uniquely determined up to
$*$-isomorphism, \cite[Th.~4.1]{Ham}. Injective C*-algebras are monotone
complete, and so is $I(A)$, i.e.~any bounded increasingly directed net of
self-adjoint elements of $I(A)$ admits a least upper bound in $I(A)$
(cf.~\cite{ChEf:77}).
In particular, $I(A)$ is unital. Note that the injective envelope of a
non-unital C*-algebra equals the injective envelope of its unitization.
A C*-algebra $B$ containing $A$ as a C*-subalgebra equals $I(A)$ if and only
if $B$ is injective and the restriction of any completely isometric linear map
$\phi: B \to I(A)$ to $A \subseteq B$ is completely isometric again,
cf.~\cite[Prop.~4.7]{Ham}. For more information on injective envelopes of
C*-algebras in various categories we refer to \cite{Ham,BP,FP}.

Injective commutative C*-algebras $A=C(X)$ can be characterized easily:
$A=C(X)$ is injective if and only if $X$ is stonean (\cite[Th.~7.3]{Dix}).
This follows also from \cite[Th.~25.5.1]{Sem} and the fact that bounded linear
maps between C*-algebras are positive whenever their norm equals the norm of
their evaluation at the identity of the C*-algebra.
What is more, by \cite[Th.~6.3]{Ham2} the equality $x = \sup \{ a \in A_{sa} :
a \leq x \}$ holds for every self-adjoint element $x$ of the injective envelope
$I(A)$ of a commutative C*-algebra $A$. Since for $x \in I(A)_{sa}^+$ and
$A$ commutative the inequality $a \leq x$ for some $a \in A_{sa}$ implies
$0 \leq (a \vee 0) \leq x$ the set $\{ a \in A_{sa}^+ : a \leq x \}$ is
non-trivial for any positive $x \in I(A)$, and the supremum of this set equals
exactly $x$. By \cite[Lemma 1.7]{Ham3} the latter condition is equivalent
to the assertion that every non-zero projection of the injective envelope $I(A)$
of a commutative C*-algebra $A$ mayorizes some non-zero positive element of
$A$. Analyzing this inner order structure of an injective envelope $I(A)$ of a
unital commutative C*-algebra $A=C(X) \,$ H.~Gonshor identified $I(A)$ with the
C*-algebra $B(X)$ of all bounded Borel functions on $X$ modulo the ideal of
all functions supported on meager subsets of $X$, cf.~\cite[Th.~1]{Go2}. The
background of this coincidence is that every remainder class of $B(X)$
contains a lower semi-continuous on $X$ function (\cite[Lemma 7.5]{Dix}).

\bigskip
The second main construction in the present paper is also essentially an
algebraic one - the local multiplier algebra $M_{loc}(A)$ of a given
C*-algebra $A$.
A norm-closed one-(resp., two-)sided ideal $I$ of $A$ is said to be {\it
essential} if its intersection with any other one-(resp., two-)sided ideal of
$A$ is non-trivial. For commutative C*-algebras $A=C_0(X)$ an ideal $I$ is
essential if and only if there exists a closed meager subset $S_I \subseteq X$
such that $I=\{ f \in C_0(X) : f \,\,\, {\rm is} \,\, {\rm continuous} \,\,
{\rm on} \,\,\, X \setminus S_I \}$, where $X$ is the locally compact Hausdorff
space corresponding to $A$ by I.~M.~Gelfand's theorem.
If for an arbitrary C*-algebra $A$ and any pair of essential two-sided ideals
$J \subseteq I \subseteq A$ their multiplier C*-algebras $M(J) \supseteq M(I)
\supseteq M(A)$ are considered then the induced partial order of canonical
embeddings of those multiplier algebras $\{ M(I) : I \,\,\, {\rm essential}
\,\, \mbox{\rm two-sided} \,\, {\rm in} \,\, \, A \}$ gives rise to the direct
limit of this partially ordered set:
  \[
    Q_b(A) = {\rm alg} \,\lim_\rightarrow \{ M(I) : I \,\,\, {\rm essential} \,\,
    \mbox{\rm two-sided} \,\, {\rm ideal} \,\, {\rm in} \,\,\, A \}
  \]
Performing the construction in the set $B(H)$ of all bounded linear operators
on a Hilbert space $H$ wherein $A$ is faithfully represented, the norm-closure
$M_{loc}(A)$ of $Q_b(A)$ becomes a unital C*-algebra. It is said to be the
{\it local multiplier algebra of} $A$, cf.~\cite{Ell,Ped78,AM}.
Note that the local multiplier algebra of a non-unital C*-algebra equals the
local multiplier algebra of its unitization.
If $A=C(X)$ is commutative then $M_{loc}(C(X))$ can be identified with the set
of all continuous functions over the inverse limit $\lim_\leftarrow \{ \beta U
: U \,\,\, {\rm open} \,\, {\rm dense} \,\, {\rm in} \,\,\, X \}$ which is a
compact Hausdorff space, \cite[Th.~1]{AM}. Moreover, it is a commutative
AW*-algebra by \cite[Th.~1]{AM}, and so it is monotone complete by
\cite[Prop.~2.2]{Az}. For more information on local multiplier algebras we
refer to \cite{AM,AM2,So1,So2}.

%%%%%%%%%%%%%%%%%%%%%%%%%%%%%%%%%%%%%%%%%%%%%%%%%%%%%%%%%%%%%%%%%%%%%%%%%%%%%%%

\section{Results}

The connection between the local multiplier algebra and the injective envelope
of C*-algebras is described by the following theorem:

\medskip
{\bf Theorem 1:}          %  \label{included}
  {\sl Let $A$ be a (unital) C*-algebra, $I(A)$ be its injective envelope and
  $M_{loc}(A)$ be its local multiplier algebra. Then $M_{loc}(A)$ can be
  canonically identified with a C*-subalgebra of $I(A)$, i.e.~the canonical
  embedding $A \hookrightarrow I(A)$ extends to an isometric algebraic
  monomorphism $M_{loc}(A) \hookrightarrow I(A)$. 
  If $A$ is commutative then $M_{loc}(A)$ coincides with $I(A)$.  }

\medskip
{\sl Proof:}
Without loss of generality we assume $A$ to be unital.
The C*-algebra $M_{loc}(A)$ is an operator $A$-$A$ bimodule by its construction,
and $A$ is a C*-subalgebra of $M_{loc}(A)$. Consider the canonical embedding
of $A$ into $I(A)$. By the injectivity of $I(A)$ and by G.~Wittstock's
extension theorems for completely bounded $A$-$A$ bimodule maps
(\cite{Wittstock}) this canonical
embedding $A \subseteq I(A)$ extends to a completely contractive $A$-$A$ bimodule
map $\psi: M_{loc}(A) \to I(A)$ with the same complete boundedness norm one,
see \cite[Th.~3.1]{Wittstock}, \cite{Ar:69}. Let us show that this extension
is actually an isometric algebraic embedding of $M_{loc}(A)$ into $I(A)$.

Consider an essential two-sided ideal $I \subseteq A$ and its multiplier
algebra $M(I) \subseteq M_{loc}(A)$. Every element $x \in M(I)$ is representable
as the strict limit of a net of elements $\{ x_\alpha \}$ of
$I$, i.e.
  \[
    \lim_{\alpha} \| xy-x_\alpha y \| =
    \lim_{\alpha} \| yx-yx_\alpha \| = 0
  \]
for every $y \in I$. Therefore, inside $I(A)$ we obtain the equalities
  \begin{eqnarray*}
    \psi(x)y & = & \psi(xy) = \lim_{\alpha} \psi(x_\alpha y)
    = \lim_{\alpha} \psi(1_A) x_\alpha y = xy   \\
    y\psi(x) & = & \psi(yx) = \lim_{\alpha} \psi(yx_\alpha)
    = \lim_{\alpha} \psi(1_A) y x_\alpha = yx
  \end{eqnarray*}
for every $y \in I$, where all limits are limits in norm. So $\psi(x)$ is the
strict limit of the net $\{ x_\alpha \} \subset I \subseteq I(A)$ with respect
to $I$ since $I(A)$ is an AW*-algebra and this strict limit is realized in a
unique way as an element of $I(A)$ by \cite[Th.]{Ped}. Consequently,
every multiplier algebra $M(I)$ of an essential two-sided ideal $I$ of $A$
is isometrically algebraically embedded by the existing map $\psi$ into
$I(A)$, and the selected embedding preserves canonical containment relations
between essential two-sided ideals of $A$ , the C*-algebra $A$ itself and
their multiplier algebras. So the map $\psi$ is actually an injective
$*$-homomorphism of C*-algebras.

\smallskip
If $A$ is commutative then $M_{loc}(A)$ is a commutative AW*-algebra by
\cite[Th.~1]{AM}. Hence, $M_{loc}(A)$ is an injective C*-algebra, and by
W.~B.~Arveson's extension theorem (\cite{Ar:69}) for completely positive maps
applied to the identity map on $A$ and by the definition of injective
envelopes and their uniqueness, the identity $I(A) \equiv M_{loc}(A)$ holds.  $\bullet$

By \cite[Cor.~1.6]{Ham3} there are C*-algebras $A$ such that the injective
envelope $I(Z(A))$ of the center $Z(A)$ of $A$ is different from the center
$Z(I(A))$ of the injective envelope $I(A)$ of $A$. We identify $I(Z(A))$
with the center $Z(M_{loc}(A))$ of $M_{loc}(A) \subseteq I(A)$ that belongs
to the center of $I(A)$ whenever $Z(I) \circ I$ is dense in $I$ for any
essential two-sided ideal $I \subseteq A$.

\medskip
{\bf Theorem 2:}    %  \label{center}
   {\sl Let $A$ be a C*-algebra with a center $Z(A)$ such that $Z(A) \circ A$ is
   norm-dense in $A$. Then the center of the multiplier C*-algebra of $A$
   coincides with the multiplier C*-algebra of the center of $A$,
   i.e.~$Z(M(A)) \equiv M(Z(A))$.
   \newline
   If, moreover, $Z(I) \circ I$ is norm-dense in $I$ for any essential
   two-sided ideal $I$ of $A$, then the center of the local multiplier
   C*-algebra of $A$ coincides with the local multiplier C*-algebra of the
   center of $A$, i.e.~$Z(M_{loc}(A)) \equiv M_{loc}(Z(A))$. Moreover,
   $M_{loc}(Z(A)) \equiv I(Z(A))$. }

\medskip
{\sl Proof:}
At the beginning we show that the center of the multiplier C*-algebra of any
C*-algebra $B$ coincides with the multiplier C*-algebra of its center,
i.e.~$Z(M(B)) \equiv M(Z(B))$, as soon as $Z(B) \circ B$ is norm-dense in $B$.
The inclusion $M(Z(B)) \subseteq Z(M(B))$ is
immediate. Indeed, $M(Z(B))$ consists of all strict limits of strictly converging
norm-bounded nets of $Z(B)$. The multiplier C*-algebra $M(B)$ is a Banach
$Z(M(B))$-$Z(M(B))$ bimodule such that both $Z(M(B)) \circ M(B)$ and $M(B)
\circ Z(M(B))$ are norm-dense in $M(B)$. Hence, \cite[Th.~4.1]{Ped97} applies
saying that every element $x \in M(B)$ can be represented as $x=z_1y_1$ and as
$x = y_2z_2$ for certain positive
$z_1,z_2 \in Z(M(B))$ with $\|z\| \leq 1$ and certain $y_1,y_2 \in M(B)$.
Let $\{ z_\alpha \}$ be a norm-bounded strictly converging net of $Z(A)$.
Then for any $x \in M(B)$ both the nets $\{ z_\alpha x \}$ and $\{ x z_\alpha
\}$ converge in norm inside the Banach $Z(M(B))$-$Z(M(B))$-module $M(B)$
since $z_\alpha x = (z_\alpha z_1) y_1$, $x z_\alpha = y_2 (z_2 z_\alpha)$
and the coefficients at the right sides are converging in norm. Therefore,
every strictly converging norm-bounded net of $Z(B)$ preserves strict
convergence if considered as a norm-bounded net of $M(B)$. So its strict limit
belongs to $M(B)$ and, moreover, to $Z(M(B))$ since strict convergence
preserves commutation relations in the limits.

Conversely, fix $t \in Z(M(B))$ and an approximate
identity $\{s_\alpha \}$ of $Z(B)$. We have $(ts_\alpha)x = txs_\alpha =
x(ts_\alpha)= x(s_\alpha t)$ for any $x \in B$. So $ts_\alpha = s_\alpha t
\in Z(B)$ for any index $\alpha$. What is more, $(ts_\alpha)u = t(s_\alpha u)$
for any $u \in Z(B)$. The right side converges in norm, so does the left side
showing the strict convergence of the net $\{ ts_\alpha \} \in Z(B)$. Since
$Z(B) \circ B$ is norm-dense in $B$ the net $\{ ts_\alpha \} \in Z(B)$
converges strictly to the element $t \in M(Z(B))$.

The C*-algebra $B$ can be replaced by any essential two-sided norm-closed
ideal $I$ of the given C*-algebra $A$ by assumption. The additional condition claimed
above ensures that the idea of the first part of the present proof can be
applied again. If two essential ideals $I,J$ of $A$
fulfil $I \subseteq J$ then their centers satisfy $Z(I) \subseteq Z(J)$.
This can be seen switching to the respective multiplier algebras $M(J)
\subseteq M(I)$. We derive
\begin{eqnarray*}
   Z(\, {\rm alg} \,\, \lim_\to \,\, M(I)\, ) & = &
        {\rm alg} \,\, \lim_\to \,\, Z(M(I)) \\
        & = & {\rm alg} \,\, \lim_\to \,\, M(Z(I))
\end{eqnarray*}
for the algebraic direct limit of the partially ordered net of multiplier
C*-algebras of these essential ideals $I \subseteq A$. If $J$ is an essential
norm-closed ideal of $Z(A)$ then the norm-closure of $A J A$ is an essential
two-sided norm-closed ideal of $A$ with center $J$. So at the right side the
multiplier C*-algebra of any essential two-sided norm-closed ideal $J$ of $Z(A)$
appears, and we obtain
\[
    Z(\, {\rm alg} \,\, \lim_\to \,\, M(I)\, ) =
          {\rm alg} \,\, \lim_\to \,\, M(J)
\]
for $I \subseteq A$ and $J \subseteq Z(A)$ as described. Forming the norm-closure
of both the sides of the latter equality we get $M_{loc}(Z(A))$ at the
right side and $Z(M_{loc}(A))$ at the left by \cite[Th.~1, Lemma 1, Cor.~1]{AM}.
This shows the first set identity claimed. The other identity of C*-algebras
is a simple consequence of Theorem 1 applied to the center $Z(A)$
of the C*-algebra $A$.  $\bullet$

Theorem 2 is surely not true for $A=K(H)$ and $M_{loc}(A) =
M(A) = B(H)$, where $H=l_2$ is the standard separable Hilbert space.
Similar problems appear if one of the essential two-sided ideals $I$ of a
certain C*-algebra $A$ has a center $Z(I)$ for which $Z(I) \circ I$ does not
cover one or more block-diagonal components of $I$. In connection with this
phenomenon it would be interesting to know whether the additional conditions
on $A$ can be weakened for Theorem 2 to hold, or not. Unfortunately,
we are not able to formulate any conjecture on the diversity of the appearing
situations at present. Nevertheless, Theorem 2 shows the way to a
new short proof of \cite[Cor.~2]{AM} in a particular situation:

\medskip
{\bf Corollary 3:}  % \label{cor3}
   {\sl Let $A$ be a C*-algebra with the property that $Z(I) \circ I$ is norm-dense
   in $I$ for any essential two-sided ideal $I$ of $A$. Then we have the set
   identity $Z(M_{loc}(M_{loc}(A))) \equiv Z(M_{loc}(A))$ of the respective
   centers, i.e.~taking higher order local multiplier C*-algebras does not
   change the center any more. }

\medskip
{\sl Proof:}
By the Theorems 1 and 2 and by the local multiplier
C*-algebra properties of injective C*-algebras we conclude
\begin{eqnarray*}
   Z(M_{loc}(M_{loc}(A))) & \equiv & M_{loc}(Z(M_{loc}(A)))
   \equiv M_{loc}(I(Z(A))) \\
   & \equiv & I(Z(A)) \equiv Z(M_{loc}(A))  \, ,
\end{eqnarray*}
cf.~\cite[Th.]{Ped}.   $\bullet$

\smallskip
The observation formulated at Theorem 1 shows the way to a
criterion on the coincidence of the (local) multiplier algebra of a given
C*-algebra $A$ with the injective envelope of $A$ in terms of a Wittstock
type extension theorem. Recall the definition of {\it operator $B$-$C$
bimodules} ${}_B {\mathcal M}_C$
over unital C*-algebras $B$ and $C$, \cite{Wittstock,BMP}. Operator spaces
$\mathcal M$ can be characterized as norm-closed subspaces of C*-algebras.
Then operator $B$-$C$ bimodules are operator spaces $\mathcal M$ which are
$B$-$C$ bimodules such that the trilinear module pairing $B \times \mathcal M
\times C \to \mathcal M \, , \,\, (b, x, c) \to bxc$ is completely contractive
in the sense of E.~Christensen and A.~Sinclair \cite{C-S}, i.e.~it extends to
be a completely contractive linear map on the Haagerup tensor product $B
\otimes_h \mathcal M \otimes_h C$. We may assume that all module actions are
unital.

\medskip
{\bf Proposition 4:}
  {\sl Let $A$ be a C*-algebra, $B,C$ be unital C*-subalgebras of $M(A)$ (resp., of
  $M_{loc}(A)$). The following two conditions are equivalent:

  \noindent
  \begin{tabular}{cp{13cm}}
  (i) &  For any operator $B$-$C$ bimodule $\mathcal M$, any operator $B$-$C$
       subbimodule $\mathcal N$ and any completely bounded $B$-$C$ bilinear map
       $\phi: \mathcal N \to M(A)$ (resp., $\phi: \mathcal N \to M_{loc}(A)$)
       there exists a completely bounded $B$-$C$ bilinear map
       $\psi: \mathcal M \to M(A)$ (resp.,          
       $\psi: \mathcal M \to M_{loc}(A)$)
       such that $\psi |_{\mathcal N} = \phi$ and $\| \psi \|_{c.b.} =
       \| \phi \|_{c.b.}$.   \\
  (ii)&  The C*-algebras $M(A)$ and $I(A)$ (resp., $M_{loc}(A)$ and $I(A)$)
       are $*$-isomorphic.
  \end{tabular}
  }

\medskip
{\sl Proof:}
Condition (ii) implies condition (i) by G.~Wittstock's extension theorem, see
\cite{Wittstock,MN}. Conversely, consider the canonical $*$-monomorphism of
$M(A)$ (resp., $M_{loc}(A)$) into $I(A)$ existing by Theorem 1.
Both $M(A)$ (resp., $M_{loc}(A)$) and $I(A)$ are operator $B$-$C$ bimodules
if considered in a faithful $*$-representation of $I(A)$.
By \cite[Th.~2.1]{FP} the identity map of $A$ onto its canonical copy inside
$I(A)$ admits a unique completely bounded extension to $I(A)$ with the same
complete boundedness norm in case we require one-sided $A$-linearity: the
identity map on $I(A)$. If condition (i) is supposed to hold then the
canonical copy of $M(A)$ (resp., $M_{loc}(A)$) inside $I(A)$ has to coincide
with $I(A)$. This implies (ii). $\bullet$

\smallskip
We do not know whether the process of taking higher order local multiplier
algebras of a C*-algebra $A$ in general ever stabilizes (or at least converges), or
not. However, it takes place in a fixed monotone complete C*-algebra, the
injective envelope $I(A)$ of $A$, and more precisely in the monotone
closure $\overline{A}$ of $A$ in $I(A)$ that is a uniquely determined
monotone complete C*-algebra by \cite{Ham2}. $\overline{A}$ is said to be
the {\it regular monotone completion of $A$}. We remark that C*-subalgebras
$A$ of injective C*-algebras $B$ might not admit an embedding of their
injective envelopes $I(A)$ as a C*-subalgebra of $B$ that extends the given
embedding of $A$ into $B$. An example was given by M.~Hamana,
\cite[Rem.~3.9]{Ham2}.

\medskip
{\bf Theorem 5:}  % \label{prop}
  {\sl Let $A$ be a C*-algebra, $M(A)$ be its multiplier algebra and $M_{loc}(A)$
  be its local multiplier C*-algebra.
  Then the injective envelopes $I(A)$, $I(M(A))$ and $I(M_{loc}(A))$ of these
  three C*-algebras, respectively, coincide. Every higher order local multiplier
  algebra of $A$ is a C*-subalgebra of the regular monotone completion
  $\overline{A}$ of $A$ inside $I(A)$. }

\medskip
{\sl Proof:}
Since $M(A)$ and $M_{loc}(A)$ are $*$-isomorphically embedded into $I(A)$
extending the canonical embedding of $A$ into $I(A)$, the C*-algebra $I(A)$
serves as an injective extension of both $M(A)$ and $M_{loc}(A)$. However,
the identity map on $M(A)$ or $M_{loc}(A)$, respectively, admits a unique
extension to a completely positive map on $I(A)$ of the same complete
boundedness norm one because $A \subseteq M_{loc}(A) \subseteq I(A)$ by
construction and $I(A)$ is the injective envelope of $A$ by definition,
cf.~\cite[Th.~2.1]{FP}. So $I(A)$ serves as the injective envelope of $M(A)$
and $M_{loc}(A)$, too.

By \cite[Th.]{Ped} the multiplier C*-algebra of any C*-subalgebra $B$ of the
regular monotone completion $\overline{A}$ of $A$ inside $I(A)$ can be
realized as a C*-subalgebra of $\overline{A}$. So the process of taking the
higher order local multiplier C*-algebras of $A$ can be always carried out inside
$\overline{A}$. $\bullet$

\smallskip
D.~W.~B.~Somerset pointed out to us some non-standard examples of C*-algebras
$A$ beyond the set of commutative C*-algebras for which $M_{loc}(A) \equiv
M_{loc}(M_{loc}(A)))$, \cite{So2}. However, the class of C*-algebras that
coincide with their own local multiplier C*-algebra is very heterogeneous from
the point of view of existing classifications. So it also contains all AW*-algebras
(like von Neumann algebras) together with all unital simple C*-algebras
(like $C_r^*(F_2)$, irrational $A_\theta$, ${\mathcal O}_n$). For example,
while in AW*-algebras the linear span of the set of projections is norm-dense
the C*-algebra $C_r^*(F_2)$ does not contain any non-trivial projection.
Nevertheless, every C*-algebra of this class possesses a commutative
AW*-algebra as its center, cf.~\cite{Dav,Ped,AM,AM2}.

\medskip
{\bf Problem 6:}
   {\sl Whether any C*-algebra $A$ possess the property $M_{loc}(A) \equiv
   M_{loc}(M_{loc}(A)))$, and if not, does the process of consecutively taking
   local multiplier C*-algebras stabilize? Are there intrinsic characterizations
   of C*-algebras $A$ for which $A \equiv M_{loc}(A)$?  }

%The next statement generalizes a fact first observed by D.~W.~B.~Somerset
%in particular situations in \cite{So2}. It establishes the set of all local
%multiplier C*-algebras as a class of C*-algebras of special interest. 
%
%\begin{proposition}  \label{locked}
%   For any C*-algebra $A$ we have $M_{loc}(A) \equiv M_{loc}(M_{loc}(A))$.
%\end{proposition}
%
%\begin{proof}
%Consider an essential two-sided norm-closed ideal $J$ of $M_{loc}(A)$. Our
%aim is to show that its multiplier C*-algebra can be realized as a C*-subalgebra
%of $M_{loc}(A)$ itself. This would show the assertion.
%
%The intersection $I = J \cap A$ is a non-zero two-sided
%norm-closed ideal of $A$. Let us show that it is essential for $A$. Let $I_0$
%be another two-sided norm-closed ideal of $A$ with $I_0 \cap I = \{ 0 \}$.
%Then the central carrier projections $q,q_0$ of $I,I_0$, respectively, in the
%bidual von Neumann algebra $A^{**}$ of $A$ are orthogonal to one another.
%Since there is a canonical $*$-monomorphism of $M(I)$ into $A^{**}$ for any
%norm-closed ideal $I \subseteq A$ that preserves the commutant we obtain
%a canonical $*$-monomorphism $M_{loc}(A) \hookrightarrow A^{**}$ that extends
%the canonical embedding of $A$ into $A^{**}$, and $J \subseteq qA^{**}$.
%However, the norm-closure of the set $(M_{loc}(A) I_0 M_{loc}(A))$ is a
%two-sided norm-closed ideal of $M_{loc}(A)$ contained in $q_0A^{**}$. So
%$J$ is essential in $M_{loc}(A)$ if and only if $I_0$ is trivial in $A$,
%and the ideal $I$ is essential in $A$.
%
%To derive the assertion above we note that $M(J) \subseteq M(I) \subseteq
%M_{loc}(A)$.
%\end{proof}

\medskip
For non-commutative C*-algebras $A$ with non-injective multiplier C*-algebra
$M(A)$ the injective envelope may be strictly larger than $M_{loc}(A)$ since
some of the essential one-sided ideals may have larger algebras of left
multipliers enjoying an isometric algebraic embedding into $I(A)$ by
G.~Wittstock's extension theorem for completely bounded one-sided module maps
\cite[Th.~4.1]{Wittstock}, the properties of the left-strict (resp.,
right-strict) topology on the algebra of left multipliers $LM(A)$ (resp., right
multipliers $RM(A)$) and by the monotone completeness of $I(A)$.

However, even if $M(I) \equiv LM(I)$ for any two-sided ideal $I$ of a certain
C*-algebra $A$ (as valid for any AW*-algebra - see \cite[Th.]{Ped}) the
coincidence of $M_{loc}(A)$ and $I(A)$ may fail. Any non-injective von Neumann
algebra $A$ serves as an example.

\medskip
{\bf Corollary 7:}
   {\sl Let $A$ be an AW*-algebra and $I \subseteq A$ be an essential two-sided
   norm-closed ideal of $A$. For any C*-algebra $B$ with $I \subseteq B
   \subseteq A$ we obtain $M_{loc}(B) \equiv A$. The C*-algebra $A$ is
   injective if and only if $M_{loc}(B) \equiv I(B) \equiv A$.  }

\medskip
{\sl Proof:}
Since $I \subseteq B$ and $I$ is an essential two-sided norm-closed ideal of
$A$ the C*-algebra $I$ is an essential two-sided norm-closed ideal of $B$ and
$M(I) \equiv A$. So $A \subseteq M_{loc}(B)$. At the other side the multipliers
of any C*-subalgebra $J$ of an AW*-algebra $A$ can be found inside $A$ by
\cite[Th.]{Ped}. In other words, $A$ always contains a faithful $*$-representation
of $M(J)$ that extends the given embedding of $J$ into $A$. Hence, $M_{loc}(B)
\subseteq A$. The coincidence of both these C*-algebras follows.
$M_{loc}(B)=A$ is injective if and only if $A$ is injective, otherwise
$I(B) \equiv I(M_{loc}(B)) \equiv I(A) \supset A$ by Theorem 5.
$\bullet$

\smallskip
To indicate some concrete examples, we can consider the pair of C*-algebras
$I=K(H)$ and $A=B(H)$ of all compact linear operators on an infinite-dimensional
Hilbert space $H$ and of all bounded linear operators on $H$, the pair of an
injective ${\rm II}_\infty$ von Neumann factor $A$ and its unique non-trivial
two-sided norm-closed ideal $I \subset A$, or the pair of the C*-algebra
$I=C_0((0,1])$ of all continuous functions on the unit interval that vanish at
zero and the Dixmier algebra $A=D([0,1])$ which is defined to be the set of all
bounded Borel functions on $[0,1]$ modulo the ideal of all functions supported
on meager subsets of $[0,1]$.

Finally, we arrive at the following open questions which are surely difficult
to answer:

\medskip
{\bf Problem 8:}
   {\sl Characterize the C*-algebras $A$ for which the local
   multiplier C*-algebra $M_{loc}(A)$ of $A$ coincides with the injective
   envelope $I(A)$ of $A$, or at least with the regular monotone completion
   $\overline{A}$ of $A$ in $I(A)$.  }

\medskip
Combining Theorem 2 and Corollary 3 with
\cite[Cor.~1.6]{Ham3} we may obtain necessary conditions to the center of
certain C*-algebras $A$ to satisfy $M_{loc}(A) \equiv \overline{A}$ and
$M_{loc}(A) \equiv I(A)$.
Investigating AW*-algebras $A$, for example, the local multiplier algebra
$M_{loc}(A)$ of $A$ coincides with $A$ itself by \cite[Th.]{Ped}. However,
$A$ coincides with its regular monotone completion $\overline{A}$ if and only
if $A$ is monotone complete. So we arrive at a long standing open problem of
C*-theory dating back to the work of I.~Kaplansky in 1951 (\cite{Kapl}): Are
all AW*-algebras monotone complete, or do there exist counterexamples?

\smallskip
Consider the universal $*$-representation of $A$ in $B(H)$. By the
Arveson-Wittstock theorem \cite{Ar:69,Wittstock} the canonical embedding
of $A$ into its bidual von Neumann algebra $A^{**}$ realized as the bicommutant
$A''$ of $A$ in $B(H)$ extends to a completely isometric linear map $\psi$ of
$I(A)$ into $B(H)$. The image $\psi(I(A)) \subseteq B(H)$ is an operator
system. The domain of $\psi$ contains $M_{loc}(A)$ as a C*-algebra on which
$\psi$ acts as an algebraic $*$-monomorphism by Theorem 1,
canonical identifications provided.

\medskip
{\bf Problem 9:}
   {\sl Is $M_{loc}(A)$ the largest C*-subalgebra of $I(A)$ on which $\psi$ acts
   as an algebraic $*$-monomorphism ? }

\bigskip
{\bf Acknowledgements:} The author is grateful to K.~R.~Davidson, M.~Ma\-thieu,
Zhong-Jin Ruan and D.~W.~B.~Somerset for valuable comments during discussions.
He thanks D.~P.~Blecher and V.~I.~Paulsen for their warm hospitality and
fruitful scientific collaboration during a one year stay at the University of
Houston in 1998.

%%%%%%%%%%%%%%%%%%%%%%%%%%%

%%%%%%%%%%%%%%%%%%%%%%%%%%%

\end{document}